\documentclass[12pt,a4paper,reqno]{amsart}
\usepackage{amsthm}
\usepackage{amsmath}
\usepackage{amsfonts}
\usepackage{amssymb}
\usepackage[left=2.8cm, right=2.8cm, top=3.0cm, bottom=3.0cm]{geometry}
\usepackage{color}
\usepackage[arrow, matrix, curve]{xy}
\usepackage[colorlinks=true,linkcolor=red,anchorcolor=blue,citecolor=green]{hyperref}

\def\be{\begin{equation}}
\def\ee{\end{equation}}
\def\bea{\begin{eqnarray}}
\def\eea{\end{eqnarray}}
\def\bes{\begin{eqnarray*}}
\def\ees{\end{eqnarray*}}

\def\nn{\nonumber}
\def\lb{\label}
\def\bs{\setminus}

\def\T{{\mathcal{T}}}
\def\H{{\mathcal{H}}}
\def\R{{\bf R}}
\def\C{{\bf C}}
\def\Z{{\bf Z}}

\def\N{{\bf N}}
\def\U{{\bf U}}

\def\Q{{\bf Q}}
\def\T{{\bf T}}

\def\Sg{{\Sigma}}

\def\aa{{\alpha}}
\def\bb{{\beta}}
\def\ga{{\gamma}}

\def\th{{\theta}}

\def\om{{\omega}}
\def\Om{{\Omega}}
\def\ep{{\epsilon}}
\def\lm{{\lambda}}

\def\sg{{\sigma}}
\def\dm{{\diamond}}
\def\Sg{{\Sigma}}
\def\vf{{\varphi}}

\def\<{{\langle}}
\def\>{{\rangle}}
\def\T{{\mathcal{T}}}

\def\P{{\mathcal{P}}}

\def\Sp{{\rm Sp}}

\def\hb{\vrule height0.18cm width0.14cm $\,$}

\title[Closed characteristics on compact star-shaped hypersurfaces]{Four closed characteristics on compact star-shaped hypersurfaces in $\R^{8}$}

\author[Huagui Duan]{Huagui Duan}
\address{(Huagui Duan) School of Mathematical Sciences and LPMC, Nankai University, Tianjin 300071, The People's Republic of China}
\email{duanhg@nankai.edu.cn.}

\author[Dong Xie]{Dong Xie}
\address{(Dong Xie) School of Mathematical Sciences, Nankai University,
Tianjin 300071, The People's Republic of China}
\email{1120200022@mail.nankai.edu.cn.}

\date{\today}

\begin{document}
\maketitle

\begin{abstract}
{\it In this paper, we proved that for every non-degenerate $C^3$ compact star-shaped hypersurface $\Sigma$ in $\R^{8}$ which carries no prime closed characteristic of Maslov-type index $-1$, there exist at least four prime closed characteristics on $\Sigma$.}
\end{abstract}

\renewcommand{\theequation}{\thesection.\arabic{equation}}
\renewcommand{\thefigure}{\thesection.\arabic{figure}}

\baselineskip 18pt

\setcounter{figure}{0}
\setcounter{equation}{0}
\section{Introduction and main results}

Let $\Sigma$ be a $C^3$ compact hypersurface in $\R^{2n}$ strictly star-shaped with respect to the origin, i.e.,
the tangent hyperplane at any $x\in\Sigma$ does not intersect the origin. We denote the set of all such hypersurfaces
by $\H_{st}(2n)$, and denote by $\H_{con}(2n)$ the subset of $\H_{st}(2n)$ which consists of all strictly convex
hypersurfaces. We consider closed characteristics $(\tau, y)$ on $\Sigma$, which are solutions of the following
problem
\be
\left\{\begin{array}{ll}\dot{y}=JN_\Sigma(y), \\
 y(\tau)=y(0), \lb{1.1}
\end{array}\right.
\ee
where $J=\left(\begin{array}{ll}0 &-I_n\\
 I_n & 0\end{array}\right)$, $I_n$ is the identity matrix in $\R^n$, $\tau>0$, $N_\Sigma(y)$ is the outward
normal vector of $\Sigma$ at $y$ normalized by the condition $N_\Sigma(y)\cdot y=1$. Here $a\cdot b$ denotes
the standard inner product of $a, b\in\R^{2n}$. A closed characteristic $(\tau, y)$ is {\it prime}, if $\tau$
is the minimal period of $y$. Two closed characteristics $(\tau, y)$ and $(\sigma, z)$ are {\it geometrically
distinct}, if $y(\R)\not= z(\R)$. We denote by $\T(\Sigma)$ the set of geometrically distinct
closed characteristics $(\tau, y)$ on $\Sigma\in\mathcal{H}_{st}(2n)$. A closed characteristic
$(\tau,y)$ is {\it non-degenerate} if $1$ is a Floquet multiplier of $y$ of precisely algebraic multiplicity
$2$; {\it hyperbolic} if $1$ is a double Floquet multiplier of it and all the other Floquet multipliers
are not on ${\bf U}=\{z\in {\bf C}\mid |z|=1\}$, i.e., the unit circle in the complex plane; {\it elliptic}
if all the Floquet multipliers of $y$ are on ${\bf U}$. We call a $\Sigma\in \mathcal{H}_{st}(2n)$ {\it
non-degenerate} if all the closed characteristics on $\Sigma$, together with all of their iterations, are
non-degenerate.

There is a long-standing conjecture on the number of closed characteristics on compact convex hypersurfaces in
$\R^{2n}$:
\be \,^{\#}\T(\Sg)\ge n, \qquad \forall \; \Sg\in\H_{con}(2n). \lb{1.2}\ee

In 1978, P. Rabinowitz in \cite{Rab1} proved
$^\#\T(\Sg)\ge 1$ for any $\Sg\in\H_{st}(2n)$, A. Weinstein in \cite{Wei1}
proved $^\#\T(\Sg)\ge 1$ for any $\Sg\in\H_{con}(2n)$ independently.
When $n\ge 2$, in 1987-1988, I. Ekeland-L. Lassoued, I. Ekeland-H. Hofer, and A. Szulkin (cf. \cite{EkL1},
\cite{EkH1}, \cite{Szu1}) proved $\,^{\#}\T(\Sg)\ge 2$ for any $\Sg\in\H_{con}(2n)$.

In \cite{LoZ} of 2002, Y. Long and C. Zhu further proved
\bea\;^{\#}\T(\Sg)\ge\left[\frac{n}{2}\right]+1, \qquad \forall\, \Sg\in \H_{con}(2n). \nn\eea
In particular, if all the prime closed characteristics on $\Sg$ are non-degenerate, then $\,^{\#}\T(\Sg)\ge n$
(cf. Theorem 1.1 and Corollary 1.1 of \cite{LoZ}). In \cite{WHL} of 2007, W. Wang, X. Hu and Y. Long proved
$\,^{\#}\T(\Sg)\ge 3$ for every $\Sg\in\H_{con}(6)$. In \cite{Wan1} of 2016, W. Wang proved
$\,^{\#}\T(\Sg)\ge \left[\frac{n+1}{2}\right]+1$ for every $\Sg\in\H_{con}(2n)$. In \cite{Wan2} of 2016, W.
Wang proved $\,^{\#}\T(\Sg)\ge 4$ for every $\Sg\in\H_{con}(8)$.

Note that every contact form supporting the standard contact structure on $M=S^{2n-1}$ arises from embeddings of $M$
into $\R^{2n}$ as a strictly star-shaped hypersurface enclosing the origin, it is conjectured that in fact the
conjecture (\ref{1.2}) holds for any $\Sg\in\H_{st}(2n)$ too. For the star-shaped case, \cite{Gir1} of 1984 and
\cite{BLMR} of 1985 show that $\;^{\#}\T(\Sg)\ge n$ for $\Sg\in\H_{st}(2n)$ under some pinching conditions.
In \cite{Vit2} of 1989, C. Viterbo proved a generic existence result for infinitely many closed
characteristics on star-shaped hypersurfaces in ${\bf R}^{4n}$. In \cite{HuL} of 2002, X. Hu and Y. Long proved that
$\;^{\#}\T(\Sg)\ge 2$ for any non-degenerate $\Sg\in \H_{st}(2n)$. In \cite{HWZ2} of 2003, H. Hofer, K. Wysocki, and
E. Zehnder proved that $\,^{\#}\T(\Sg)=2$ or $\infty$ holds for every non-degenerate $\Sg\in\H_{st}(4)$ provided that
the stable and unstable manifolds of every hyperbolic closed orbit on $\Sg$ intersect transversally. This
condition was removed by D. Cristofaro-Gardiner, M. Hutchings and D. Pomerleano \cite{CGHP} of 2019.
Recently, D. Cristofaro-Gardiner, U. Hryniewicz, M. Hutchings and H. Liu in \cite{CGHHL} proved that $\,^{\#}\T(\Sg)=2$ or $\infty$ holds for every $\Sg\in\H_{st}(4)$. In \cite{CGH1} of 2016, D. Cristofaro-Gardiner and M. Hutchings proved that $\;^{\#}\T(\Sg)\ge 2$ for any contact three
manifold $\Sg$. Various proofs of this result for a star-shaped hypersurface can be found in \cite{GHHM}, \cite{LLo1}
and \cite{GiG1}.

In \cite{GuK} of 2016, J. Gutt and J. Kang proved $\,^{\#}\T(\Sg)\ge n$ for every non-degenerate
$\Sg\in \H_{st}(2n)$ if every closed characteristic on $\Sg$ possesses the Conley-Zehnder index of at least $(n-1)$. In
\cite{DLLW} of 2018, H. Duan, H. Liu, Y. Long and W. Wang proved $\,^{\#}\T(\Sg)\ge n$ for every index perfect
non-degenerate $\Sg\in \H_{st}(2n)$, and there exist at least $n$ (or $(n-1)$) non-hyperbolic closed characteristics
when $n$ is even (or odd). Here $\Sg\in \H_{st}(2n)$ is {\it index perfect} if it carries only finitely many
geometrically distinct prime closed characteristics, and every prime closed characteristic $(\tau,y)$ on $\Sigma$
possesses positive mean index and whose Maslov-type index $i(y, m)$ of its $m$-th iterate satisfies $i(y, m)\not= -1$
when $n$ is even, and $i(y, m)\not\in \{-2,-1,0\}$ when $n$ is odd for all $m\in\N$. Later V. Ginzburg, B. G\"{u}rel
and L. Macarini in \cite{GGM} extended these results to prequantization bundles which covered many known results for
closed geodesics on Finsler manifolds.

This paper continues to study the above conjecture on non-degenerate star-shaped $\Sigma$ in $\R^{8}$. Motivated by some ideas from \cite{DLLW} and \cite{DLLW2}, by removing the assumption of the positive mean index and relaxing the index restriction $i(y, m)\neq -1$ with $m\in\N$ to $i(y,1)\neq -1$ for each prime $y$ in the index perfect condition, we can still prove the existence of four prime closed characteristics on such $\Sigma$.

\medskip

{\bf Theorem 1.1.} {\it Let $\Sigma$ be a $ C^{3} $ compact non-degenerate star-shaped hypersurface in $\R^{8}$. If for every prime closed characteristic $(\tau,y)$ on $\Sigma$, the Maslov-type index of $(\tau,y)$ satisfies $i(y,1)\not= -1$, then there exist at least $4$ geometrically distinct closed characteristics. }

\medskip

In order to prove Theorem 1.1, we assume that such $\Sigma\subseteq\R^8$ possesses only finitely many prime closed characteristics $\{(\tau_k,y_k)\}_{k=1}^q$ with $i(y_{k},1)\neq -1$.

On one hand, it can be shown that the assumption of $i(y_{k},1)\neq -1$ with any $1\le k\le q$ implies that $i(y_{k}^{m})\neq -5 $ or $i(y_{k},m)\neq -1 $ for $m\in\N$ (see Lemma 3.1).

On the other hand, by using the index iteration theory and some precise analysis, we can show that if further the mean index $\hat{i}(c_k)=0$, then its Viterbo index must satisfy $i(y_{k}^{m})\in\{-6, -4\}$ (see Claim 1), and then, through following main ideas from the proof of Theorem 1.6 in \cite{DLLW2} (essentially originated from \cite{Vit2}), we conclude that this is impossible. Hence there always hold $\hat{i}(c_k)\neq 0$ for any $1\le k\le q$ (see Lemma 3.2).

In a word, under the assumption in Theorem 1.1, we can show that all conditions in Theorem 1.2 in \cite{DLLW2} can be satisfied. Therefore we obtain the existence of at least four geometrically distinct closed characteristics on such $\Sigma$.

\medskip

In this paper, let $\N$, $\N_0$, $\Z$, $\Q$, $\R$, $\C$ and $\R^+$ denote the sets of natural integers,
non-negative integers, integers, rational numbers, real numbers, complex numbers and positive real
numbers respectively. We define the functions
\bea [a]&=&\max{\{k\in {\bf Z}\mid k\leq a\}},\qquad \{a\}=a-[a],\nn\\
E(a)&=&\min{\{k\in{\bf Z}\mid k\geq a\}},\ \varphi(a)=E(a)-[a].\lb{1.3}\eea
Denote by $a\cdot b$ and $|a|$ the standard inner product and norm in $\R^{2n}$. Denote by $\langle\cdot,\cdot\rangle$ and $\|\cdot\|$
the standard $L^2$ inner product and $L^2$ norm. For an $S^1$-space $X$, we denote by
$X_{S^1}$ the homotopy quotient of $X$ by $S^1$, i.e., $X_{S^1}=S^\infty\times_{S^1}X$,
where $S^\infty$ is the unit sphere in an infinite dimensional complex Hilbert space.
We use $\Q$ coefficients for all homological and cohomological modules. By $t\to a^+$, we
mean $t>a$ and $t\to a$.

\setcounter{figure}{0}
\setcounter{equation}{0}
\section{Variational structure for closed characteristics and index iteration theory}
\subsection{Variational structure for closed characteristics}

Fix a $\Sg\in\H_{st}(2n)$, we assume the following condition on $\T(\Sg)$:

\medskip

{\bf (F)} {\it There exist only finitely many geometrically distinct prime closed characteristics
$\{(\tau_j, y_j)\}_{1\le j\le k}$ on $\Sigma$.}

\medskip

Let $j: \R^{2n}\rightarrow\R$ be the gauge function of $\Sigma$, i.e., $j(\lambda x)=\lambda$ for $x\in\Sigma$
and $\lambda\ge0$, then $j\in C^3(\R^{2n}\bs\{0\}, \R)\cap C^0(\R^{2n}, \R)$ and $\Sigma=j^{-1}(1)$.
Let $\hat{\tau}=\inf_{1\leq j\leq k}{\tau_j}$ and $T$ be a fixed positive constant. Then following \cite{Vit1}
and Section 2 of \cite{LLW}, for any $a>\frac{\hat{\tau}}{T}$, we can construct a function
$\varphi_a\in C^{\infty}({\bf R}, {\bf R}^+)$ which has $0$ as its unique critical point in $[0, +\infty)$. Moreover,
$\frac{\varphi_a^{\prime}(t)}{t}$ is strictly decreasing for $t>0$ together with $\varphi_a(0)=0=\varphi_a^{\prime}(0)$
and $\varphi_a^{\prime\prime}(0)=1=\lim_{t\rightarrow 0^+}\frac{\varphi_a^{\prime}(t)}{t}$. The precise definition of
$\varphi_a$ and the dependence of $\varphi_a$ on $a$ are given in Lemma 2.2 and Remark 2.3 of \cite{LLW} respectively.
As in \cite{LLW}, we define a Hamiltonian function
$H_a\in C^{3}({\bf R}^{2n} \setminus\{0\},{\bf R})\cap C^{1}({\bf R}^{2n},{\bf R})$ satisfying $H_a(x)=a\vf_a(j(x))$
on $U_A=\{x\mid a\vf_a(j(x))\leq A\}$ for some large $A$, and $H_a(x)=\frac{1}{2}\ep_a|x|^2$ outside some even larger
ball with $\ep_a>0$ small enough such that outside $U_A$ both $\nabla H_a(x)\not= 0$ and $H_a^{\prime\prime}(x)<\ep_a$
hold.

We consider the following fixed period problem
\be \dot{x}(t)=JH_a^\prime(x(t)),\quad x(0)=x(T). \lb{2.1}\ee
Then solutions of (\ref{2.1}) are $x\equiv 0$ and $x=\rho y(\tau t/T)$ with
$\frac{\vf_a^\prime(\rho)}{\rho}=\frac{\tau}{aT}$, where $(\tau, y)$ is a solution of (\ref{1.1}). In particular,
non-zero solutions of (\ref{2.1}) are in one to one correspondence with solutions of (\ref{1.1}) with period
$\tau<aT$.

For any $a>\frac{\hat{\tau}}{T}$, we can choose some large constant $K=K(a)$ such that
\be H_{a,K}(x) = H_a(x)+\frac{1}{2}K|x|^2 \lb{2.2}\ee
is a strictly convex function, that is,
\be (\nabla H_{a, K}(x)-\nabla H_{a, K}(y), x-y) \geq \frac{\ep}{2}|x-y|^2, \quad \forall x,y\in\R^n, \lb{2.3}\ee
for some positive $\ep$. Let $H_{a,K}^*$ be the Fenchel dual of $H_{a,K}$
defined by
$$ H_{a,K}^\ast (y) = \sup\{x\cdot y-H_{a,K}(x)\;|\; x\in \R^{2n}\}. $$
The dual action functional on $X=W^{1, 2}({\bf R}/{T {\bf Z}}, {\bf R}^{2n})$ is defined by
\be F_{a,K}(x) = \int_0^T{\left[\frac{1}{2}(J\dot{x}-K x,x)+H_{a,K}^*(-J\dot{x}+K x)\right]dt}. \lb{2.4}\ee
Then $F_{a,K}\in C^{1,1}(X, \R)$ and for $KT\not\in 2\pi{\bf Z}$, $F_{a,K}$ satisfies the
Palais-Smale condition and $x$ is a critical point of $F_{a, K}$ if and only if it is a solution of (\ref{2.1}). Moreover,
$F_{a, K}(x_a)<0$ and it is independent of $K$ for every critical point $x_a\neq 0$ of $F_{a, K}$.

When $KT\notin 2\pi{\bf Z}$, the map $x\mapsto -J\dot{x}+Kx$ is a isomorphism between
$X=W^{1, 2}({\bf R}/({T {\bf Z}}); {\bf R}^{2n})$ and $E=L^{2}({\bf R}/(T {\bf Z}),{\bf R}^{2n})$. We denote its inverse
by $M_K$ and
\be \Psi_{a,K}(u)=\int_0^T{\left[-\frac{1}{2}(M_{K}u, u)+H_{a,K}^*(u)\right]dt}, \qquad \forall\,u\in E. \lb{2.5}\ee
Then $x\in X$ is a critical point of $F_{a,K}$ if and only if $u=-J\dot{x}+Kx$ is a critical point of $\Psi_{a, K}$.

Suppose $u$ is a nonzero critical point of $\Psi_{a, K}$.
Then the formal Hessian of $\Psi_{a, K}$ at $u$ is defined by
\be Q_{a,K}(v)=\int_0^T(-M_K v\cdot v+H_{a,K}^{*\prime\prime}(u)v\cdot v)dt, \lb{2.6}\ee
which defines an orthogonal splitting $E=E_-\oplus E_0\oplus E_+$ of $E$ into negative, zero and positive subspaces.
The index and nullity of $u$ are defined by $i_K(u)=\dim E_-$ and $\nu_K(u)=\dim E_0$ respectively.
Similarly, we define the index and nullity of $x=M_Ku$ for $F_{a, K}$, we denote them by $i_K(x)$ and
$\nu_K(x)$. Then we have
\be i_K(u)=i_K(x),\quad \nu_K(u)=\nu_K(x), \lb{2.7}\ee
which follow from the definitions (\ref{2.4}) and (\ref{2.5}). The following important formula was proved in
Lemma 6.4 of \cite{Vit2}:
\be i_K(x) = 2n([KT/{2\pi}]+1)+i^v(x) \equiv d(K)+i^v(x), \lb{2.8}\ee
where the Viterbo index $i^v(x)$ does not depend on $K$, but only on $H_a$.

By the proof of Proposition 2 of \cite{Vit1}, we have that $v\in E$ belongs to the null space of $Q_{a, K}$
if and only if $z=M_K v$ is a solution of the linearized system
\be \dot{z}(t) = JH_a''(x(t))z(t). \lb{2.9}\ee
Thus the nullity in (\ref{2.7}) is independent of $K$, denoted by $\nu^v(x)\equiv \nu_K(u)= \nu_K(x)$.

By Proposition 2.11 of \cite{LLW}, the index $i^v(x)$ and nullity $\nu^v(x)$ coincide with those defined for
the Hamiltonian $H(x)=j(x)^\alpha$ for all $x\in\R^{2n}$ and some $\aa\in (1,2)$. Especially
$1\le \nu^v(x)\le 2n-1$ always holds.

For every closed characteristic $(\tau, y)$ on $\Sigma$, let $aT>\tau$ and choose $\vf_a$ as above.
Determine $\rho$ uniquely by $\frac{\vf_a'(\rho)}{\rho}=\frac{\tau}{aT}$. Let $x=\rho y(\frac{\tau t}{T})$.
Then we define the index $i(\tau,y)$ and nullity $\nu(\tau,y)$ of $(\tau,y)$ by
\bea
i(\tau,y)=i^v(x), \qquad \nu(\tau,y)=\nu^v(x).\lb{2.9'}
\eea
Then the mean index of $(\tau,y)$ is defined by
\bea \hat i(\tau,y) = \lim_{m\rightarrow\infty}\frac{i(m\tau,y)}{m}. \nn\eea
Note that by Proposition 2.11 of \cite{LLW}, the index and nullity are well defined and are independent of the
choice of $a$. For a closed characteristic $(\tau,y)$ on $\Sigma$, we simply denote by $y^m\equiv(m\tau,y)$
the m-th iteration of $y$ for $m\in\N$.

\subsection{The index iteration theory for symplectic paths}

In \cite{Lon2} of 1999, Y. Long established the basic normal form
decomposition of symplectic matrices. Based on this result he
further established the precise iteration formulae of indices of
symplectic paths in \cite{Lon3} of 2000.

As in \cite{Lon3}, denote by
\bea
N_1(\lm, b) &=& \left(\begin{array}{cc}\lm & b\\
 0 & \lm \end{array}\right), \qquad {\rm for\;}\lm=\pm 1, \; b\in\R, \lb{3.1}\\
D(\lm) &=& \left(\begin{array}{cc}\lm & 0\\
 0 & \lm^{-1} \end{array}\right), \qquad {\rm for\;}\lm\in\R\bs\{0, \pm 1\}, \lb{3.2}\\
R(\th) &=& \left(\begin{array}{cc}\cos\th & -\sin\th \\
 \sin\th & \cos\th \end{array}\right), \qquad {\rm for\;}\th\in (0,\pi)\cup (\pi,2\pi), \lb{3.3}\\
N_2(e^{\th\sqrt{-1}}, B) &=& \left(\begin{array}{cc} R(\th) & B \\
 0 & R(\th) \end{array} \right), \qquad {\rm for\;}\th\in (0,\pi)\cup (\pi,2\pi)\;\; {\rm and}\; \nn\\
 && \quad B=\left(\begin{array}{cc} b_1 & b_2\\
 b_3 & b_4 \end{array}\right)\; {\rm with}\; b_j\in\R, \;\;
 {\rm and}\;\; b_2\not= b_3. \lb{3.4}\eea
Here $N_2(e^{\th\sqrt{-1}}, B)$ is non-trivial if $(b_2-b_3)\sin\theta<0$, and trivial
if $(b_2-b_3)\sin\theta>0$.

As in \cite{Lon3}, the $\diamond$-sum (direct sum) of any two real matrices is defined by
$$ \left(\begin{array}{cc}A_1 & B_1\\ C_1 & D_1 \end{array}\right)_{2i\times 2i}\diamond
 \left(\begin{array}{cc}A_2 & B_2\\ C_2 & D_2 \end{array}\right)_{2j\times 2j}
=\left(\begin{array}{cccc}A_1 & 0 & B_1 & 0 \\
 0 & A_2 & 0& B_2\\
 C_1 & 0 & D_1 & 0 \\
 0 & C_2 & 0 & D_2\end{array}\right). $$

For every $M\in\Sp(2n)$, the homotopy set $\Omega(M)$ of $M$ in $\Sp(2n)$ is defined by
$$ \Om(M)=\{N\in\Sp(2n)\,|\,\sg(N)\cap\U=\sg(M)\cap\U\equiv\Gamma,
 \;\nu_{\om}(N)=\nu_{\om}(M),\, \forall\om\in\Gamma\}, $$
where $\sg(M)$ denotes the spectrum of $M$,
$\nu_{\om}(M)\equiv\dim_{\C}\ker_{\C}(M-\om I)$ for $\om\in\U$.
The component $\Om^0(M)$ of $P$ in $\Sp(2n)$ is defined by
the path connected component of $\Om(M)$ containing $M$.

\medskip

For every $\ga\in\mathcal{P}_\tau(2n)\equiv\{\ga\in C([0,\tau],Sp(2n))\ |\ \ga(0)=I_{2n}\}$, we extend
$\ga(t)$ to $t\in [0,m\tau]$ for every $m\in\N$ by
\bea
\ga^m(t)=\ga(t-j\tau)\ga(\tau)^j \qquad \forall\;j\tau\le t\le (j+1)\tau \;\;
{\rm and}\;\;j=0, 1, \ldots, m-1,\lb{3.8}
\eea
as in p.114 of \cite{Lon2}. As in \cite{LoZ} and \cite{Lon4}, we denote the Maslov-type indices of
$\ga^m$ by $(i(\ga,m),\nu(\ga,m))$. Note that these differences of these index and nullity from those defined
in (\ref{2.9'}) are given in Theorem 2.2 below.

The following is the precise index iteration formulae for symplectic paths, which is due to Y. Long (cf. Theorems
1.2 and 1.3 of \cite{Lon3} and Theorems 8.2.1 and 8.3.1 in Chapter 8 of \cite{Lon4}).

\medskip

{\bf Theorem 2.1.} {\it Let $\ga\in\P_{\tau}(2n)$. Then there exists a path $f\in C([0,1],\Omega^0(\gamma(\tau))$
	such that $f(0)=\gamma(\tau)$ and
	\bea f(1)
	&=& N_1(1,1)^{\diamond p_-} \diamond I_{2p_0}\diamond N_1(1,-1)^{\diamond p_+}
	\diamond N_1(-1,1)^{\diamond q_-} \diamond (-I_{2q_0})\diamond N_1(-1,-1)^{\diamond q_+}\nn\\
	&& \diamond R(\theta_1)\diamond\cdots\diamond R(\theta_r)\diamond N_2(\omega_1, u_1)\diamond\cdots
	\diamond N_2(\omega_{r_*}, u_{r_*}) \nn\\
	&& \diamond N_2(\lm_1, v_1)\diamond\cdots\diamond N_2(\lm_{r_0}, v_{r_0})\diamond M_s \nn\eea
	where $N_2(\omega_j, u_j)$s are non-trivial and $N_2(\lm_j, v_j)$s are trivial basic normal forms; $\sigma(M_s)\cap U=\emptyset$
	with $M_s=D(2)^{\diamond s}$ with an integer $s\ge 0$ or $M_s=D(-2)\diamond D(2)^{\diamond (s-1)}$ with an integer $s\ge 1$;
 $\omega_j=e^{\sqrt{-1}\alpha_j}$, $\lambda_j=e^{\sqrt{-1}\beta_j}$; $\theta_j$, $\alpha_j$, $\beta_j\in (0, \pi)\cup (\pi, 2\pi)$; $p_-$, $p_0$, $p_+$, $q_-$, $q_0$, $q_+$, $r$, $r_*$ and $r_0$ are non-negative integers and satisfy the equality
 \bea
 p_{-}+p_{0}+p_{+}+q_{-}+q_{0}+q_{+}+r+2r_{*}+2r_{0}+s=n\lb{3.9}.
 \eea
 These integers and real numbers are uniquely determined by $\gamma(\tau)$. Then we have
	\bea i(\gamma, m)
	&=& m(i(\gamma,1)+p_-+p_0-r)+2\sum_{j=1}^r E\left(\frac{m\theta_j}{2\pi}\right)-r-p_--p_0\nn\\
	& & -\frac{1+(-1)^m}{2}(q_0+q_+)+2\left(\sum_{j=1}^{r_*}\varphi\left(\frac{m\alpha_j}{2\pi}\right)-r_*\right), \lb{3.9'}
	\eea
	\bea \nu(\gamma, m)
	&=& \nu(\gamma,1)+\frac{1+(-1)^m}{2}(q_-+2q_0+q_+)+2(r+r_*+r_0)\nn\\
	& & -2\left(\sum_{j=1}^{r}\varphi\left(\frac{m\theta_j}{2\pi}\right)+\sum_{j=1}^{r_*}\varphi\left(\frac{m\alpha_j}{2\pi}\right)
	+\sum_{j=1}^{r_0}\varphi\left(\frac{m\beta_j}{2\pi}\right)\right),\lb{3.9''}\\
	\hat{i}(\gamma, 1) &=& \lim_{m\to +\infty}\frac{i(\gamma,m)}{m}=i(\gamma, 1)+p_-+p_0-r+\sum_{j=1}^r\frac{\theta_j}{\pi},\lb{3.9'''}
	\eea

	We have $i(\gamma, 1)$ is odd if $f(1)=N_1(1, 1)$, $I_2$, $N_1(-1, 1)$, $-I_2$,
	$N_1(-1, -1)$ and $R(\theta)$; $i(\gamma, 1)$ is even if $f(1)=N_1(1, -1)$ and $N_2(\omega, b)$; $i(\gamma,1)$
	can be any integer if $\sigma(f(1))\cap\U = \emptyset$.}

\medskip

{\bf Theorem 2.2.} (cf. Theorem 2.1 of \cite{HuL} and Theorem 6.1 of \cite{LLo2}) {\it Suppose $\Sg\in \H_{st}(2n)$ and
	$(\tau,y)\in \T(\Sigma)$. Then we have
	\be i(y^m)\equiv i(m\tau,y)=i(y, m)-n,\quad \nu(y^m)\equiv\nu(m\tau, y)=\nu(y, m),
	\lb{3.32}\ee
	where $m\in\N$, $i(y^m)$ and $\nu(y^m)$ are the index and nullity of $(m\tau,y)$ defined in Section 2.1, $i(y, m)$ and $\nu(y, m)$
	are the Maslov-type index and nullity of $(m\tau,y)$ (cf. Section 5.4 of \cite{Lon3}). In particular, we have
	$\hat{i}(\tau,y)=\hat{i}(y,1)$, where $\hat{i}(\tau, y)$ is given in Section 2.1, $\hat{i}(y,1)$
	is the mean Maslov-type index. Next we denote it simply by $\hat{i}(y)$.}

\setcounter{figure}{0}
\setcounter{equation}{0}
\section{Proof of Theorem 1.1}

In this section, let $\Sigma\in \mathcal{H}_{st}(8)$ be a non-degenerate compact star-shaped hypersurface in $\R^8$. Next we assume that such $\Sigma$ carries only finitely many prime closed characteristics $\{(\tau_k,y_k)\}_{k=1}^q$ with $i(y_{k},1)\neq -1$.

Let $P_{\Sg} = \{m\tau_k\;|\;1\le k\le q, m\in\N\}$ be the period set of all closed characteristics.
Denote by $\ga_k\equiv \ga_{y_k}$ the associated symplectic path of $(\tau_k,y_k)$ for $1\le k\le q$.
Then by Lemma 3.3 of \cite{HuL} and Lemma 3.2 of \cite{Lon1}, there exists $P_k\in Sp(8)$ and $U_k\in Sp(6)$
such that
\bea M_k\equiv\ga_k(\tau_k)=P_k^{-1}(N_1(1,1)\dm U_k)P_k,\qquad \forall\ 1\le k\le q,\lb{5.1}\eea
where, because $\Sigma$ is non-degenerate, every $U_k$ is isotopic in $\Omega^0(U_k)$ to a matrix of the following form by Theorem 2.1 (cf. Theorem 4.7 of \cite{Wan2})
\bea
&& R(\th_{k,1})\,\dm\,\cdots\,\dm\,R(\th_{k,r_{k}})\,\dm\,D(\pm 2)^{\dm s_{k}} \nn\\
&& \qquad\dm\,N_2(e^{\aa_{k,1}\sqrt{-1}},A_{k,1})\,\dm\,\cdots\,\dm\,N_2(e^{\aa_{k,r_{k,\ast}}\sqrt{-1}},A_{k,r_{k,\ast}})\nn\\
&& \qquad \dm\,N_2(e^{\bb_{k,1}\sqrt{-1}},B_{k,1})\,\dm\,\cdots\,\dm\,N_2(e^{\bb_{k,r_{k,0}}\sqrt{-1}},B_{k,r_{k,0}}), \nn\eea
where $\frac{\th_{k,j}}{2\pi}\in[0,1]\setminus\Q$ for $1\le j\le r_{k}$; $\frac{\aa_{k,j}}{2\pi}\in[0,1]\setminus\Q$ for $1\le j\le r_{k,\ast}$;
$\frac{\bb_{k,j}}{2\pi}\in[0,1]\setminus\Q$ for $1\le j\le r_{k,0}$ and
\be
r_{k}+ s_{k} +2r_{k,\ast} + 2r_{k,0} = 3,\qquad \forall\ 1\le k\le q. \lb{5.2}
\ee
Hence by (\ref{5.1}), Theorem 2.2 and the precise index iteration formulae for symplectic paths in Theorem 2.1 , for $1\le k\le q$, noticing that here $n=4$ and $p_-=1$, we have
\bea
i(y_k^m)
&=& m(i(y_k)+n+1-r_{k})+2\sum_{j=1}^{r_{k}} \left[\frac{m\theta_{k,j}}{2\pi}\right]+r_{k}-1-n\nn\\
&=& m(i(y_k)+5-r_{k})+2\sum_{j=1}^{r_{k}} \left[\frac{m\theta_{k,j}}{2\pi}\right]+r_{k}-5,\quad\forall\ m\ge 1, \lb{5.3}
\eea
where in (\ref{5.3}), we have used $E(a)=[a]+1$ for $a\in \R\bs\Z$. Thus
\be
\hat{i}(y_k)=i(y_k)+5-r_{k}+\sum_{j=1}^{r_{k}}\frac{\theta_{k,j}}{\pi}, \quad\forall\ 1\le k\le q. \lb{5.4}
\ee

An iterate $(m\tau, y)$ of a prime closed characteristic $(\tau,y)$ on $\Sigma$ with $m\in\N$ is called
{\it good}, if its Maslov-type index has the same parity as that of $(\tau,y)$.

\medskip

{\bf Lemma 3.1.} {\it For any $ 1\le k\le q $, the Viterbo index of each good $ m $-th iterate $ y_{k}^{m} $ of $ y_{k} $ with $ m \in \N $ satisfies $ i(y_{k}^{m})\neq -5 $.}

\medskip

{\bf Proof.} Suppose there exists $1 \leq k_0 \leq q$ and an integer $m_0 \geq 2$ such that the iterate $y_{k_0}^{m_0}$ is good and $i(y_{k_0}^{m_0}) = -5$. According to Theorem 2.2 and the assumption $i(y_{k_0},1)\neq -1$, we obtain that $i(y_{k_0}) \in 2\Z + 1 \setminus \{-5\}$. Hence, the estimate holds
\bea
|i(y_{k_{0}})+5|\ge 2.\lb{A1}
\eea
Using (\ref{5.3}), we obtain
\bea
m_{0}(i(y_{k_{0}})+5-r_{k_{0}})+2\sum_{j=1}^{r_{k_{0}}} \left[\frac{m_{0}\theta_{k_{0},j}}{2\pi}\right]+r_{k_{0}}=0.\lb{A2}
\eea
From this we obtain
\bea
i(y_{k_{0}})+5=r_{k_{0}}-\frac{1}{m_{0}}\left(2\sum_{j=1}^{r_{k_{0}}} \left[\frac{m_{0}\theta_{k_{0},j}}{2\pi}\right]+r_{k_{0}}\right). \lb{A3}
\eea
Combining (\ref{A1}) and (\ref{A3}), it yields
\bea
\left|(m_{0}-1)r_{k_{0}}-2\sum_{j=1}^{r_{k_{0}}}\left[\frac{m_{0}\theta_{k_{0},j}}{2\pi}\right]\right| \ge 2m_{0}.\lb{A4}
\eea
Note that $0 \leq \left[\frac{m_0\theta_{k_0,j}}{2\pi} \right] \leq m_0-1$ by the definition (\ref{1.3}) and $\frac{\theta_{k_0,j}}{2\pi}\in(0,1)\setminus \Q$ , thus
\bea
\left|(m_{0}-1)r_{k_{0}}-2\sum_{j=1}^{r_{k_{0}}}\left[\frac{m_{0}\theta_{k_{0},j}}{2\pi}\right]\right| \le (m_{0}-1)r_{k_{0}}.\lb{A5}
\eea
Combining (\ref{A4}), (\ref{A5}), and (\ref{5.2}), we conclude that
\bea
r_{k_{0}}=3.\lb{A6}
\eea

According to Theorems 2.3 and 2.4, (\ref{A6}) implies that $i(y_{k_0}) \in 2\Z$. This leads to a contradiction.\hfill\hb

\medskip

{\bf Lemma 3.2.} {\it For any $ 1\le k\le q $, there holds $\hat{i}(y_k)\neq 0$.}

\medskip

{\bf Proof.} We assume that there exist an integer $ q_{0}\in[1,q] $ such that $\hat{i}(y_k)=0$ for any $1\le k\le q_{0}$. Then we have the following claim.

\medskip

{\bf Claim 1.} {\it For $1\le k\le q_{0}$ and for any good $ m $-th iterate $ y_{k}^{m} $ with some $ m\in\N $, the following index estimate of $y_{k}^{m}$ holds
\bea
	i(y_{k}^{m})\in\{-6, -4\}.
\eea}

{\bf Proof.} By (\ref{5.4}) we have
\be
\hat{i}(y_k)=i(y_k)+5-r_{k}+\sum_{j=1}^{r_{k}}\frac{\theta_{k,j}}{\pi}=0,\quad \forall\ 1\le k\le q_0. \lb{5.5}
\ee

Note that $\frac{\theta_{k,j}}{\pi}\notin \Q$ and $0\leq r_{k}\leq 3$ by (\ref{5.2}), and then it yields $r_k\neq 1$ by (\ref{5.5}). Therefore we have \bea r_{k}\in\{0,2,3\}.\lb{3.13a}\eea

By (\ref{5.3}) and (\ref{5.5}), for any $ m\in\N $ we have
\bea
i(y_{k}^{m}) &=& m\left(-\sum_{j=1}^{r_{k}}\frac{\theta_{k,j}}{\pi}\right)+2\sum_{j=1}^{r_{k}} \left[\frac{m\theta_{k,j}}{2\pi}\right]+r_{k}-5\nn\\
&=& -2\sum_{j=1}^{r_{k}}\left\{\frac{m\theta_{k,j}}{2\pi}\right\}+r_{k}-5.\lb{3.50a}
\eea
When $ r_{k}\in\{2,3\}$, noting that $\frac{\theta_{k,j}}{\pi}\notin \Q$, by (\ref{3.50a}) there holds
\bea \sum_{j=1}^{r_{k}}\{\frac{m\theta_{k,j}}{2\pi}\}\in [1,r_{k}-1]\cap\N.\lb{3.51}\eea

More precisely, if $ r_{k}=2 $, by (\ref{3.50a}) and (\ref{3.51}) it yields $i(y_{k}^{m})=-5$, which contradicts to Lemma 3.1.

If $ r_{k}=3 $, it follow from (\ref{3.50a}) and (\ref{3.51}) that
\bea
i(y_{k}^{m}) = -2\sum_{j=1}^{r_{k}}\left\{\frac{m\theta_{k,j}}{2\pi}\right\}-2\in\{-4,-6\}.\lb{3.50}
\eea

If $r_k=0$, by (\ref{3.50a}) it yields $i(y_{k}^{m})=-5$, which contradicts to Lemma 3.1 again.

Thus Claim 1 is proved.
	
\medskip

Based on Claim 1, we will follow main ideas from that of Theorem 1.6 in \cite{DLLW2} (essentially originated from \cite{Vit2}) to complete the proof of Lemma 3.2. For the convenience of readers, we only outline some necessary details.

\medskip

{\bf Step 1.} On one hand, by Claim 1, there holds $i(y_{k}^{m})\in\{-4, -6\}$ for $1\le k\le q_{0}$ and any good $ m $-th iterate $ y_{k}^{m} $ with some $ m\in\N $. On the other hand, note that $\hat{i}(y_k)>0$ (respectively $<0$) implies $i(y_k^m)\to +\infty$ (respectively $-\infty$) as
$m\to +\infty$. Thus good iterates $y_k^m$ of every $y_k$ for $q_0+1\le k\le q$ have indices satisfying $i(y_k^m)\notin \{-3,-5,-7\} $ for any large enough $m\in\N$. Therefore for large enough $a$, all the good closed characteristics $y_k^m$ for
$1\le k\le q$ with period larger than $aT$, which implies that the iterate number $m$ is very large, will have their
Viterbo indices:
\bea
\left\{\begin{array}{ll}
	&\text{either (i) equal to} -4 \;\text{or} -6, \quad \text{when}\; \hat{i}(y_j)=0, \cr
	&\text{or (ii) different from} -3,\ -5\ \text{and} -7, \quad \text{when}\; \hat{i}(y_j)\not= 0.\cr
\end{array}
\right. \lb{4.11}
\eea

\medskip

{\bf Step 2.} For $a\in\R$, let $X^-(a,K)=\{x\in X\mid F_{a,K}(x)<0\}$ with $K=K(a)$ as defined in Section 2
as well as in Section 7 of \cite{Vit2}. For any large enough positive $a<a'$, we fix the same constant $K'>0$ as that in the second Step on p.639 of \cite{Vit2}
to be sufficiently large than $K$ such that the Hamiltonian $H_{t,K'}(x)$ is strictly convex for every
$t\in [a,a']$. Now let $A=X^-(a,K')$ and $A'=X^-(a',K')$.
Because the period set $P_{\Sg}= \{m\tau_j\;|\;1\le j\le 2, m\in\N\}$ is discrete, we choose the above constants $a$
and $a'$ carefully such that $aT$ and $a'T$ do not belong to $P_{\Sg}$.

Now, for chosen sufficiently large $a$ and $a'$ with $a < a'$, by (\ref{4.11}), there exists no good closed characteristic whose period lies between $aT$ and $a'T$ possessing Viterbo index $ -3 $,$ -5 $ or $ -7 $. Therefore, following the discussion on pp.78-79 of \cite{Cha}, according to the arguments in Steps 2 and 3 in the proof of Theorem 1.6 of \cite{DLLW2}, we have
\be
 H_{S^1, d(K')-3}(A',A) \;=\; H_{S^1, d(K')-5}(A',A) \;=\; H_{S^1, d(K')-7}(A',A) \;=\; 0. \lb{4.13}
\ee

\medskip

{\bf Step 3.} Now we consider the following exact sequence of the triple $(X,A',A)$
\bea
&&\cdots\longrightarrow H_{S^1, d(K')+\mu+1}(A',A)\stackrel{i_{\mu+1*}}{\longrightarrow} H_{S^1, d(K')+\mu+1}(X,A)
\stackrel{j_{\mu+1*}}{\longrightarrow}H_{S^1, d(K')+\mu+1}(X,A')\nn\\
&&\quad\qquad\stackrel{\partial_{\mu+1*}}{\longrightarrow} H_{S^1, d(K')+\mu}(A',A)\stackrel{i_{\mu*}}{\longrightarrow} H_{S^1, d(K')+\mu}(X,A)
\stackrel{j_{\mu*}}{\longrightarrow}H_{S^1, d(K')+\mu}(X,A')\nn\\
&&\qquad\qquad\stackrel{\partial_{\mu*}}{\longrightarrow} H_{S^1, d(K')+\mu-1}(A',A)\longrightarrow\cdots, \lb{4.14}\eea
where $ \mu\in\{-4,-6\} $.

Next, we consider the following homomorphisms:
\bea
H_{S^1, d(K)+\mu+1}(X,X^-(a,K)) &{\xi_1}\atop{\longrightarrow}& H_{S^1, d(K')+\mu+1}(X,X^-(a,K')), \nn\\
H_{S^1, d(K')+\mu+1}(X,X^-(a,K')) &{\xi_2}\atop{\longrightarrow}& H_{S^1, d(K')+\mu+1}(X,X^-(a',K')), \nn\\
H_{S^1, d(K)+\mu+1}(X,X^-(a,K)) &{\xi}\atop{\longrightarrow}& H_{S^1, d(K')+\mu+1}(X,X^-(a',K')), \nn\eea
where $\xi_1$ is the homomorphism given by (7.2) of \cite{Vit2}, $j_{\mu+1*}=\xi_2$ is the homomorphism
given by the line above (7.4) of \cite{Vit2}, and $\xi= \xi_2\circ\xi_1$ is precisely the homomorphism
given by (7.4) of \cite{Vit2}. Here $\xi_1$ is an isomorphism and $\xi$ is a zero homomorphism as
proved in the Steps 1 and 2 of the proof of Theorem 7.1 in \cite{Vit2} respectively. Therefore,
$j_{\mu+1*}=\xi_2$ is also a zero homomorphism.

Thus, combining (\ref{4.13}) and (\ref{4.14}), we obtain
\be H_{S^1, d(K')+\mu+1}(X,A) = \text{Ker}(j_{\mu+1*}) = \text{Im}(i_{\mu+1*}) = i_{\mu+1*}(H_{S^1, d(K')+\mu+1}(A',A)) = 0. \lb{4.15}\ee

Now we fix the above-chosen $a'>0$ and select another large enough $a''>a'$, and enlarge the constant $K'$
chosen in Step 2 such that the Hamiltonian $H_{t,K'}(x)$ is also strictly convex for every
$t\in [a',a'']$. Then, repeating the above proof with the long exact sequence of the triple
$(X,A'',A')$ instead of $(X,A',A)$ in the above arguments, with $A''=X^-(a'',K')$, we similarly obtain
\be H_{S^1, d(K')+\mu+1}(X,A') =0. \lb{4.16}\ee

Combining (\ref{4.15}) and (\ref{4.16}), (\ref{4.14}) yields
\be 0\stackrel{\partial_{\mu+1*}}{\longrightarrow} H_{S^1, d(K')+\mu}(A',A)\stackrel{i_{\mu*}}{\longrightarrow}
H_{S^1, d(K')+\mu}(X,A)\stackrel{j_{\mu*}}{\longrightarrow}H_{S^1, d(K')+\mu}(X,A')
\stackrel{\partial_{\mu*}}{\longrightarrow} 0. \lb{4.17}\ee

\medskip

{\bf Step 4.} When $a$ increases, we always meet infinitely many good closed characteristics with some Viterbo index $\mu\in\{-4,-6\}$ due to Claim 1. For the above chosen large enough $a < a'$, there exist only finitely many closed characteristics among $\{y_k^m\ |\ 1\le k\le q_{0},m\ge 1\}$ such that their periods locate between $aT$ and $a'T$. As Step 5 in the proof of Theorem 1.6 of \cite{DLLW2}, we obtain
\bea
H_{S^1, d(K')+\mu}(A',A) \not= 0, \lb{4.18}
\eea

\medskip

{\bf Step 5.} By the exactness of the sequence (\ref{4.17}) and (\ref{4.18}), we obtain
$$ H_{S^1, d(K')+\mu}(X,A)=H_{S^1, d(K')+\mu}(A',A)\bigoplus H_{S^1, d(K')+\mu}(X,A')\neq 0. $$
Then, by our choice of $a$, $a'$ and $a''$, and replacing $(X,A',A)$ by $(X,A'',A')$ in the above arguments, similarly
we obtain
\be H_{S^1, d(K')+\mu}(X,A')\neq 0. \lb{4.19}\ee

Now on one hand, if $j_{\mu*}$ in (\ref{4.17}) is a trivial homomorphism, then by the exactness of the sequence
(\ref{4.17}) it yields
$$ H_{S^1, d(K')+\mu}(X,A') = \text{Ker}(\partial_{\mu*}) = \text{Im}(j_{\mu*})=0, $$
which contradicts (\ref{4.19}). Therefore $j_{\mu*}$ in (\ref{4.17}) is a non-trivial homomorphism.

However, on the other hand, by our discussion between (\ref{4.14}) and (\ref{4.15}) using the arguments in \cite{Vit2},
$j_{\mu*}$ in (\ref{4.17}) is a zero homomorphism. This contradiction completes the proof of Lemma 3.2. \hfill\hb

\medskip

{\bf Proof of Theorem 1.1.}

\medskip

By Lemma 3.1, Lemma 3.2, and Theorem 2.2, we can conclude that for every prime closed characteristic $(\tau, y)$ on $\Sigma$, the Maslov-type index of each good $m$-th iterate $(m\tau, y)$ of $(\tau, y)$ with some $ m \in\N $ satisfies $i(y, m) \neq -1$, and every prime closed characteristic on $\Sigma$ possesses nonzero mean index. Then, according to Theorem 1.2 in \cite{DLLW2}, there exist at least 4 geometrically distinct closed characteristics. Thus, the proof of Theorem 1.1 is complete.\hfill\hb

\medskip

\section*{Acknowledgements}

H. Duan is partially supported by the National Key R\&D Program of China (Grant No. 2020YFA0713300), the National Natural Science Foundation of China (Nos. 12271268 and 12361141812), and the Fundamental Research Funds for the Central Universities. D. Xie is partially supported by the National Natural Science Foundation of China (No. 12361141812).




\bibliographystyle{abbrv}

\end{document}